\documentclass[10pt,a4paper,twoside]{article}
\usepackage{amssymb}
\usepackage{amsmath}
\usepackage{fancyhdr}
\usepackage{amsthm}
\usepackage[T1]{fontenc}

\usepackage{afterpage}  %
\theoremstyle{plain}
\newtheorem{theorem}{Theorem}[section]

\theoremstyle{definition}

\theoremstyle{remark}

\begin{document}

\begin{center}
\vspace*{2pt}
{\Large \textbf{Fiala--Agre list of single axioms}}\\[3mm]
{\Large\textbf{for Boolean groups is wrong}}\\[30pt]
{\large \textsf{\emph{Aleksandar Krape\v{z}}}}
\\[30pt]
\end{center}
\textbf{Abstract} \\
We point out mistakes in the recent result by Fiala and Agre concerning short axioms for Boolean groups. \\[5mm]
\vspace*{5mm}

In their paper \cite{FialaAgre} Nick C. Fiala and Keith M. Agre give, among other results, the list of shortest single axioms (in the language $L = \{ \cdot, e \}$) for Boolean groups (groups of exponent $2$). They explain their choice of 1323 candidate formulas, a series of tests with programs Prover9 and Mace4 they submitted their formulas to and finaly give Theorem 6.2 in which they list 173 axioms for Boolean groups and further five candidate formulas which are proved to be axioms only for finite models.

Unfortunately, their Theorem is wrong. To see this recall the following well known theorem for Boolean groups.
\begin{theorem} 
{\em An identity $S = T$ is a theorem of the theory of Boolean groups iff for every variable $x$ from $S = T$\/, the number of occurrences of $x$ in $S = T$ is even.}
\end{theorem}
Inspecting $173 + 5$ formulas given in Theorem 6.2
we see that only $8 + 3$ of them
satisfy above condition, namely:
\[
\begin{alignedat}{3}
(\text{80R4}{\uparrow})&  &\qquad ((e \cdot xy) \cdot yz) z = x  &\qquad\qquad(\text{81R1}{\downarrow})& &\qquad e((x \cdot yz)y \cdot z) = x \\
(\text{81L2}{\downarrow})&   &\qquad e(xy \cdot z) \cdot yz = x &\qquad\qquad(\text{81R2}{\downarrow})& &\qquad ((e \cdot xy) \cdot zy)z = x  \\
(\text{81L3}{\downarrow})&  &\qquad (e(x \cdot yz) \cdot y) z = x &\qquad\qquad(\text{81R3}{\downarrow})& &\qquad e((x \cdot yz)z \cdot y) = x  \\
(\text{81M2}{\downarrow})&  &\qquad e(xy \cdot zy) \cdot z = x  &\qquad\qquad(\text{81R3}{\uparrow})& &\qquad e((xy \cdot z) \cdot yz) = x
\end{alignedat}
\]
and 
\[
\begin{alignedat}{3}
(\text{81L2})& &\qquad (ex \cdot yz)y \cdot z = x  &\qquad\qquad(\text{81M2})& &\qquad (ex \cdot y)z \cdot yz = x \\
(\text{81R1})& &\qquad (ex \cdot yz)z \cdot y = x & & &  
\end{alignedat}
\]
(Since formulas were neither numbered nor tagged in the original paper, I decided to tag them by the 'coordinates' which point to their place in Theorem 6.2. The first number is the page, then a column and finaly a position
from the top or bottom of the list.)

The rest of the formulas are actually equivalent to $x = e$\/, which is an
axiom for trivial (one--element) groups.

Fiala and Agre did not give explicitly their initial list of 1323 formulas 
but it is obvious by the end result that they started with the wrong one, despite Theorem 2.3 which requires that the candidate formulas are of the form $T = x$ where $T$ is a term with one occurrence of $x$ and $e$ each and with
two occurrences of $y$ and $z$ each. Another possibility is that the final list was inconsistenly edited, but this is just guessing. 

The authors also did not give 
the code for Prover9 and Mace4 they used in the case of Boolean groups,
but just state that the search 'similar' to the search for axioms for groups of exponent $3$ 'was carried out'. Therefore we cannot check whether there are mistakes in the code too.

It is highly unlikely that such widely used and long lasting programs as Prover9 and Mace4
are wrong, but it might be possible that the local copies used by authors have some mistakes.



\vspace{1cm}
\noindent
{Mathematical Institute \\
         of the Serbian Academy of Sciences and Arts\\
         Knez Mihailova 36\\
         11001 Beograd \\
         Serbia \\
e-mail: sasa@mi.sanu.ac.rs }

\end{document}